\documentclass[reqno,12pt]{amsart}
\usepackage{amsfonts, amssymb, amsthm, amsmath, chngcntr}
\usepackage{verbatim, graphicx, multicol, epic, floatflt, comment}
\usepackage{xcolor}
\newtheorem{bigthm}{Theorem}
\newtheorem{thm}{Theorem}[section]

\newtheorem{lem}[thm]{Lemma}
\newtheorem{cor}[thm]{Corollary}

\theoremstyle{definition}

\newtheorem{defn}[thm]{Definition}

\theoremstyle{remark}
\newtheorem*{rem*}{Remark}

\newtheorem*{ack*}{Acknowledgment}

\numberwithin{equation}{section}
\counterwithout{figure}{section}

\newcommand{\eps}{\ensuremath{\epsilon} }
\newcommand{\epso}{\ensuremath{\epsilon > 0\;}}
\newcommand{\bs}{\hfill\ensuremath{_\blacksquare} }

\newcommand{\N}{\ensuremath{\mathbb N} }

\newcommand{\Ci}{\ensuremath{C^{(i)}}}

\newcommand{\Pa}{\ensuremath{\mathbb P} }

\newcommand{\mB}{\ensuremath{\mathcal B} }

\title[Independent Alpern Tower]{An Alpern tower independent of a given partition}

\begin{document}

\author[J. Campbell]{James T. Campbell} \address[J. T.
Campbell]{Dept. Math. Sci.\\Dunn Hall 373\\
  University of Memphis\\Memphis, TN 38152}
\email{jcampbll@memphis.edu}

\author[J. Collins]{Jared T. Collins}\address[J. T. Collins]{Dept. Math. and Comp. Sci.\\Freed-Hardeman University\\Henderson TN 38340}
\email{jtcollins@fhu.edu}

\author[S. Kalikow]{Steven Kalikow}\address[S. Kalikow]{Dept. Math. Sci.\\Dunn Hall 373\\
  University of Memphis\\Memphis, TN 38152}
\email{skalikow@memphis.edu}

\author[R. King]{Raena King}\address[R. King]{Dept. Math. and Comp. Sci \\Christian Bros. University\\ Memphis, TN, 38104} \email{rking2@cbu.edu}

\author[R. McCutcheon]{Randall McCutcheon} \address[R. McCutcheon]{Dept. Math. Sci.\\Dunn Hall 373\\
  University of Memphis\\Memphis, TN 38152}
\email{rmcctchn@memphis.edu}

\begin{abstract}

  Given a measure-preserving transformation $T$ of a probability space $(X, \mathcal B, \mu)$ and a finite measurable partition $P$ of $X$, we show how to construct an Alpern tower of any height whose base is independent of the partition $P$. That is, given $N \in \N$,  there exists a Rohlin tower of height $N$, with base $B$ and error set $E$, so that $B$  is independent of $P$, and $T(E) \subset B$. \bigskip

\textit{2010 Mathematics Subject Classification}: 28D05, 37M25, 60A10 \medskip

\textit{Keywords:} Rohlin tower, Alpern tower, independent sets, measure-preserving transformation, probability space.  
\end{abstract}

\date{April 30, 2015}

\maketitle

\section{Introduction and Statement of Results}
It has long been known that, given an ergodic invertible probability measure preserving system, a Rohlin tower may be constructed with base independent of a given partition of the underlying space(\cite{Roh:52}, \cite{Roh:65}). In \cite{Alp:79}, meanwhile, S. Alpern proved a `multiple' Rohlin tower theorem (see \cite{EP:97} for an easy proof) whose full statement we will not give, but which has the following corollary of interest: 

\begin{thm}\label{thm:alp}
Let $N \in \N$ and \epso be given. For any ergodic invertible measure-preserving transformation $T$ of a Lebesgue probability space $(X, \mathcal B, \mu)$, there exists a Rohlin tower of height $N$ with base $B$ and error set $E$ with $\mu(E) < \eps$, so that $T(E) \subset B$. 
\end{thm}

A {\em Rohlin tower of height} $N$ {\em with base} $B$ {\em and error set} $E$ is characterized by the collection of sets 
$\{B, TB, \dots, T^{N-1}B, E\}$ forming a partition of $X$. If in addition $T(E) \subset B$, we shall say {\bf Alpern Tower}.
It is our goal to show that for ergodic transformations on $(X, \mB, \mu)$, given a finite measurable partition $\mathbb P$ of $X$, an Alpern tower may be constructed with base $B$ independent of $\mathbb P$. Precisely:

\begin{bigthm}
\label{mainthm}
Let $(X, \mathcal B, \mu)$ be a Lebesgue probability space, and suppose $\mathbb P$ is a finite measurable partition of $X$.   For any ergodic invertible measure-preserving transformation $T$ of $X$, $N \in \N$, 
there exists a Rohlin tower of height $N$ with base $B$ and error set $E$ such that $T(E) \subset B$ and $B$ is independent of $\mathbb P$.  
\end{bigthm}

We do not specify the size of the error set; but the process of constructing our tower makes it clear that the error set may be made arbitrarily small.

\section{Proof of main result}

For the remainder of the paper, $(X, \mathcal B, \mu)$ will be a fixed Lebesgue probability space and  $T:X \to X$ will be an invertible ergodic measure-preserving transformation on $X$. All mentioned sets will be measurable and we will adopt a cavalier attitude toward null sets. In particular, ``partition'' will typically mean ``measurable partition modulo null sets''. 
\begin{defn}
By a {\em tower over B} we will mean a set $B \subset X$, called the {\em base}, and a countable partition $B = B_1 \cup B_2 \cup \cdots$, together with their images  $T^iB_j$, $0 \le i < j$, such that the family $\{T^iB_j : 0 \le i < j\}$ consists in pairwise disjoint sets. If this family partitions $X$, we will say that the tower is {\em exhaustive}. 
\end{defn}

If a tower over $B$ is exhaustive and $B = B_N \cup B_{N+1}$, we shall speak of an {\em exhaustive Alpern tower of height} $\{N, N+1\}$, as in such a case, 
$\{B, TB, \ldots, T^{N-1}B, E=T^N B_{N+1}\}$ partitions $X$ with $T(E) \subset B$. So 
we may re-phrase Theorem 1 as: \medskip

\noindent {\em {\bf Theorem \ref{mainthm}:} Let $(X, \mathcal B, \mu)$ be a Lebesgue probability space and suppose $\mathbb P$ is a finite measurable partition of $X$.   For any ergodic invertible measure-preserving transformation $T$ of $X$, $N \in \N$, one may find an exhaustive Alpern tower of height $\{N, N+1\}$ having base independent of $\mathbb P$.} \medskip

\noindent We require a lemma (and a corollary).

\begin{lem}\label{lem:m}
Let $M \in \N$ and let $\Pa = \{P_1, \dots, P_t\}$ be a partition of $X$ with $\mu(P_i) > 0$ for each $i$. There exists a set $S$ of positive measure so that if $x \in S$ with first return $n(x) = n$, say, then $|\{x, Tx, \dots, T^{n-1}x\} \cap P_i | \ge M, 1 \le i \le t$.  
\end{lem}
\proof  For almost every $x$ we may find $K(x)$ so that for each $i$ between $1$ and $t$ we have $|\{x, Tx, \dots, T^{K(x) - 1}x\} \cap P_i| \ge M$. Since almost all of $X$ is the countable union (over $k \in \N$) of $\{x: K(x) = k\}$, there exists some fixed $K$ so that the set $A = \{x: K(x) \le K\}$ has positive measure. If $C \subset A$ has very small measure ($\mu(C) < 1/K$) then the average first-return time of $x \in C$ to $C$ is $\frac{1}{\mu(C)} > K$, so we can find $S \subset C$ with $\mu(S) > 0$  so that $S, TS, \dots, T^{K - 1}S$ are pairwise disjoint. \bs

\begin{cor}\label{cor:B}
 Let $M \in \N$ and $\Pa = \{P_1, \dots, P_t\}$ be a partition of $X$ with $\mu(P_i) > 0$ for each $i$. There is a tower having base 
$S = S_{tM} \cup S_{tM+1} \cup \cdots$ where for each $x \in S_r$, $|\{x, Tx, \dots, T^{r - 1}x\} \cap P_i| \ge M$ for all $1 \le i \le t$.  
\end{cor}
\proof Let $S$, $K$  be as in Lemma \ref{lem:m} and choose any $k \ge K$. \bs \medskip

We turn now to the proof of Theorem \ref{mainthm}. Fix a partition $\Pa = \{P_1, \dots, P_t\}$, an arbitrary natural number $N$, and $\epsilon >0$. Set $m_i = \mu(P_i)$, and assume (without loss of generality) that $0 < m_1 \le m_2 \le \dots \le m_t$. Select and fix $M > \frac{3N^3t}{m_1}$. Let $S$ be as in Corollary \ref{cor:B} for this $M$; hence $S = S_{tM} \cup S_{tM+1} \cup \cdots$.  (Some $S_i$ may be empty, of course.) For each non-empty $S_R$, partition $S_R$ by $\Pa$-name of length $R$. (Recall that $x, y$ in $S_R$ have the same $\Pa$-name of length $R$ if $T^ix$ and $T^iy$ lie in the same cell of \Pa for $0 \le i < R$.) Let $C$ be the base of one of the resulting columns; hence, every $x \in C$ has the same $\Pa$-name of length $R$ (for some $R\ge tM$), and the length $R$ orbit of each $x \in C$ meets each $P_i$ at least $M$ times. 

Partition $C$ into pieces $C^{(1)}, C^{(2)}\dots$, $C^{(t)}$ whose measures will be determined later. Then partition each $C^{(i)}$ into $N$ equal measure pieces, $C^{(i)} = C^{(i)}_1 \cup C^{(i)}_2 \cup \dots \cup C^{(i)}_N$. 

Now we fix $(R,C)$ and focus our attention on the height $R$ {\em column} over a single $C^{(i)}$ and its height $R$ {\em subcolumns} over $C^{(i)}_j$, $1\leq j\leq N$. We refer to the sets $T^rC^{(i)}$, $0 \le r < R$, as {\em levels} and to the sets $T^rC^{(i)}_j$ as {\em rungs}. We are going to build a portion of $B$ by carefully selecting some rungs from the subcolumns under consideration. As we move through the various subcolumns, we need to have gaps of length $N$ or $N+1$ between selections. Now to specifics. We want 
to have our $\Ci$-selections form a ``staircase'' of height $N$ starting at level $N^2 - N$. That is, at height $(N-1)N$, the rung over $C^{(i)}_1$ is the only one selected; at height $N(N-1) + 1$, the rung over $C^{(i)}_2$ is the only one selected; etc., so that at height $N^2-1$, the rung over $C^{(i)}_N$ is the only one selected.  

This is easy to accomplish. First, we select each base rung $C^{(i)}_j$, $j = 1, 2, \dots, N$ (i.e., the rungs in the zeroth level). Over $C^{(i)}_1$, we then select $N-1$ additional rungs with gaps of length $N$; that is, we select the rungs at heights $N$, $2N, \dots, (N-1)N$. Over $C^{(i)}_2$ we select $N-2$ rungs with gap $N$, then a rung with gap $N+1$.  We continue in this fashion, choosing one less gap of length $N$ and one more of length $N+1$ in each subsequent subcolumn. In the last subcolumn (that over $C^{(i)}_N$) we are thus choosing rungs with gaps of length $N+1$ a total of $N-1$ times.  See the left side of 
Figure \ref{pic:bottom} for the case $N = 4$.

Now we perform a similar procedure moving down from the top, so as to obtain a staircase starting at height $R-(N^2-1)$. 
Note that there are either $N$ or $N-1$ unselected rungs at the top of each subcolunm. See the right side of Figure \ref{pic:bottom}. 

\begin{figure}[hbtp]

\caption[Bottom of Tower for $N = 4$]
   {Bottom, Top of Tower for $N = 4$}
\setlength{\unitlength}{.2in}

\begin{picture}(20,20)(0,0) 
\label{pic:bottom}
\put(0, 0){$C^{(i)}_1 \hspace{4mm} C^{(i)}_2 \hspace{4mm}  C^{(i)}_3 \hspace{3.3mm} C^{(i)}_4$}
\multiput(0,2)(0,1){16}{\line(1,0){1}}
\multiput(2,2)(0,1){16}{\line(1,0){1}}
\multiput(4,2)(0,1){16}{\line(1,0){1}}
\multiput(6,2)(0,1){16}{\line(1,0){1}}
\put(4,18){\vdots}
\linethickness{1mm}
\multiput(0,2)(0,4){4}{\line(1,0){1}}
\multiput(2, 2)(0,4){3}{\line(1,0){1}}
\put(2, 15){\line(1,0){1}}
\multiput(4, 2)(0,4){2}{\line(1,0){1}}
\multiput(4, 11)(0, 5){2}{\line(1, 0){1}}
\put(4, 16){\line(1,0){1}}
\multiput(6, 2)(0,5){4}{\line(1,0){1}}

\linethickness{.2mm}
\put(13, 1.5){\vdots}
\multiput(10,3)(0,1){15}{\line(1,0){1}}
\multiput(12,3)(0,1){15}{\line(1,0){1}}
\multiput(14,3)(0,1){15}{\line(1,0){1}}
\multiput(16,3)(0,1){15}{\line(1,0){1}}
\linethickness{1mm}
\multiput(10, 3)(0,5){3}{\line(1,0){1}}
\multiput(12, 4)(0,5){3}{\line(1,0){1}}
\multiput(14, 5)(0,5){2}{\line(1,0){1}}
\put(14, 14){\line(1,0){1}}
\multiput(16,6)(0,4){2}{\line(1,0){1}}
\put(16,14){\line(1,0){1}}

\end{picture}

\end{figure}
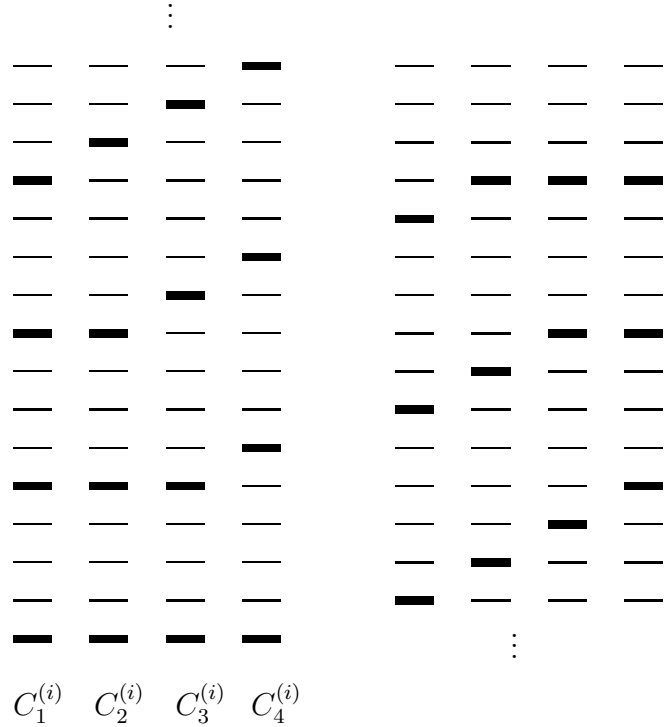 \medskip

Next, we want to select rungs through the middle of the tower so as to iterate the staircase pattern all the way up, except that we will skip certain levels (i.e. not select any of their rungs), continuing the staircase pattern where we left off with the following rung. As we want to match stride with the staircase already selected at the top, the total number of levels skipped in the middle section will be constrained to a certain residue class modulo $N$, and as we want the selected rungs to form a portion of an Alpern tower of height $\{N, N+1\}$, we cannot skip any two levels with fewer than $N$ levels between them.

Some terminology: an {\em appearance} of $P_j$ in \Ci\ is just a level of \Ci\ that is contained in $P_j$.  
A {\em selection} of $P_j$ is just a selected rung in a subcolumn of \Ci\ that is contained in $P_j$. The {\em net skips} of $P_j$ in the tower over 
\Ci\ is defined as \[ S_j(\Ci) = (\# \textnormal{ of appearances of $P_j$) } - (\# \textnormal{ of selections of $P_j$)}.\] 
For example, looking at Figure \ref{pic:bottom}, one sees that $4$ zeroth level rungs are selected. So if the zeroth level belongs to $P_j$, the zeroth level 
contribution to $S_j(\Ci)$ is $-3$ (one appearance and 4 selections).

Let $\delta = 2(N-1)(N-2)$ and choose $\gamma$ with 
\begin{displaymath}{\delta\over m_1}+N > \gamma \ge {\delta\over m_1} \quad \mbox{and} \quad (t-1)\delta + \gamma \equiv R ~(\hspace{-2.2ex}\mod N).
\end{displaymath} 
Over \Ci, we skip a quantity of ``middle'' levels belonging to each $P_j$ (for $j \neq i$) sufficient to ensure that $S_j(\Ci)=\delta$ for $j\neq i$ and $S_i(\Ci)=\gamma$. (Note that $P_j$ cannot have been skipped more than $\delta$ times in the outer rungs.) This is not delicate; one can just enact the selection greedily. That is to say, travel up the tower, beginning at level $N^2$, skipping rungs that belong to cells requiring additional skips whenever there's been no too-recent skip. Since each $P_j$ appears at least $M>{3N^3t\over m_1}$ times, and we need only $\gamma +(t-1)\delta\le {2N^2t\over m_1}$ net skips, we'll find all the skips we need.

We have not specified the relative masses of the bases of the columns $\Ci$. Set
\begin{equation}
\label{bi} b_j = \frac{\mu(P_j)(\gamma + (t-1)\delta) - \delta}{\gamma - \delta} \; 
\end{equation}
and put $\mu(\Ci)= b_i \mu(C)$, $1\leq i\leq t$. Our choice of $\gamma$ ensures that $b_i\ge 0$ for each $i$, and one easily checks that $\sum b_i=1$, so this is coherent.

Let $B_{C}$ be the union of the rungs selected from the columns over $C$ 
(this includes each of the rungs selected from each of the $N$ subcolumns over $\Ci$, $1\le i\le t$) 
and put $B=\bigcup_{C} B_{C}$. (Here $C$ runs over the bases of the columns corresponding to every $\Pa$-name of length $R$ for every $R\geq tM$.) It
is clear that $B$ forms the base of an Alpern tower of height $\{N, N+1\}$. It remains to show that $B$ is independent of $\Pa$, which we will do by constructing a set $A$, disjoint from $B$, such that both $A$ and $A\cup B$ can be shown to be independent of $\Pa$. 
 
Here is how $A$ is constructed. Consider again the tower over $\Ci$. This tower had $R$ levels and $RN$ rungs, some of which were selected for the base $B$. We now choose $\gamma +(t-1)\delta$ additional rungs for the set $A$. For each $j\neq i$, $\delta$ of these rungs should be contained in $P_j$, with the remaining $\gamma$ contained in $P_i$. (We don't worry about gaps and whatnot; just choose any such collection of rungs disjoint from the family of $B$ selections.) 
Denote the union of the these additional rungs (in all of the columns over $\Ci$, $1\leq i\leq t$) by $A_{C}$. Finally, put $A=\bigcup _{C} A_{C}$. 

That $A\cup B$ is independent of $\Pa$ is a consequence of the fact that for each $\Ci$, the number appearances of $P_j$ in the column over $\Ci$ is 
precisely the number of $B$-selections from $P_j$ plus the number of $A$-selections from $P_j$. Accordingly, the relative masses of the cells of $\Pa$ restricted to $A\cup B$ are equal to the relative frequencies of the appearances of the cells of $\Pa$ in the column over $\Ci$. Therefore, since the proportion of the column that is selected for $A\cup B$ is independent of $\Ci$ (in fact is always equal to ${1\over N}$), and since the columns over the various $\Ci$ exhaust $X$, $A\cup B$ is independent of $\Pa$ (in fact $\mu\big(P_j \cap (A\cup B)\big) = {1\over N} \mu(P_j)$, $1\leq j\leq t$).

That $A$ is independent of $\Pa$, meanwhile, is a consequence of equation (\ref{bi}). Fixing $C$ and recalling that $b_i = {\mu(\Ci)\over \mu(C)}$, 
that there were $\delta$ 
$P_j$-rungs in the column over $\Ci$ selected for $A$, $i\neq j$, and that there were $\gamma$ $P_i$-rungs in the column over $\Ci$ selected for $A$, 
the relative mass of $P_i$ among the $A$-selections in the tower over $C$ is 
\[ r_i = \frac{b_i \gamma + (1-b_i)\delta}{\gamma + (t-1)\delta} .\]
But, solving for $\mu(P_i)$ in equation (\ref{bi}), one gets that
\[ \mu(P_i) = \frac{b_i \gamma + (1-b_i)\delta}{\gamma + (t-1)\delta} \]
as well. So the intersection of $A$ with the column over $C$ is independent of $\Pa$. That this is true for every $C$ gives independence of $A$ from $\Pa$ simpliciter. \bs

\providecommand{\bysame}{\leavevmode\hbox to3em{\hrulefill}\thinspace}
\providecommand{\MR}{\relax\ifhmode\unskip\space\fi MR }
\providecommand{\MRhref}[2]{%
  \href{http://www.ams.org/mathscinet-getitem?mr=#1}{#2}
}
\providecommand{\href}[2]{#2}


\begin{thebibliography}{Roh65}

\bibitem[Alp79]{Alp:79}
Steve Alpern, \emph{Generic properties of measure preserving homeomorphisms},
  Ergodic Theory, Springer Lecture Notes in Mathematics \textbf{729} (1979),
  16--27.

\bibitem[EP97]{EP:97}
S.~J. Eigen and V.~S. Prasad, \emph{Multiple {R}okhlin tower theorem: a simple
  proof}, New York Journal of Mathematics \textbf{Proceedings of the New York
  Journal of Mathematics Conference, vol. 3A} (1997), 11--14,
  http://nyjm.albany.edu:800/j/1997/3A.11.html.

\bibitem[Roh52]{Roh:52}
V.A. Rohlin, \emph{On the fundamental ideas of measure theory}, American
  Mathematical Society Translations (1952), no.~71.

\bibitem[Roh65]{Roh:65}
\bysame, \emph{Generators in ergodic theory, {II}}, Vesnik Leningrad University
  \textbf{20} (1965), 68--72.

\end{thebibliography}
\end{document}